%% file: bdryconc.tex
\documentclass{amsart}

\input{abbrev}

\def\RR{\mathcal R}
\def\pp{\pi_1}
\def\s3{(S^3 -L)}
\def\sc{(I\times S^3 -C)}
\def\Th{\Theta}
\renewcommand{\tor}{\operatorname{Tor}}
\renewcommand{\L}{\mathcal L}
\renewcommand{\ss}{\sigma}
\renewcommand{\C}{\mathbb C}
\def\cka{\C^k_{\a}}
\def\cf{\C [F]}
\def\dl{\D_{\L}}
\begin{document}
\title[Concordance of boundary links]
{Concordance of boundary links}

\author{Jerome Levine}
\address{Department of Mathematics, 
         Brandeis University \\
         Waltham, MA 02454-9110, USA}        
\email{levine@brandeis.edu }
\urladdr{people.brandeis.edu/$\sim$levine/} 
.
 \thanks{Partially supported by an NSF grant
           and by an Israel-US BSF grant.}
\date{
\today}
\begin{abstract}
We show that the twisted signature invariants of boundary link concordance derived from unitary representations of the free group are actually ordinary link concordance invariants. We also show how the discontinuity locus of this signature function is determined by Seifert matrices of the link.
\end{abstract}

\maketitle

\section{Introduction} In \cite{L} the Atiyah-Patodi-Singer $\rho$-invariant is used to construct invariants of odd-dimensional codimension two links. The general framework is the following. For an appropriate ``universal'' group $G$ one considers the space $\RR_k  (G)$ consisting of $k$-dimensional unitary representations of $G$. If $G$ is finitely presented then $\RR_k (G)$ is a real algebraic variety. For any closed oriented odd-dimensional manifold $M$, the $\rho$-invariant defines a function $\rho_M :\RR_k (G)\to\R$, which is {\em piecewise continuous}. To obtain an invariant of a link $L$ we use $M=M_L$, the manifold obtained by doing surgery on the link components (with the $0$-framing for $1$-dimensional links). We then denote $\rho_{M_L}=\rho_L$. The group $G$ depends on the class of links considered. The concordance invariance of $\rho_L$, where the notion of concordance also depends on the choice of $G$, is expressed as follows. If $L$ is concordant to $L'$ then $\rho_L =\rho_{L'}$ in the complement of a {\em special sub-variety} of $\RR_k (G)$. For general odd-dimensional $m$-component links of dimension $n>1$, or $n=1$ with vanishing $\bar\mu$-invariants, the appropriate group $G$ is the {\em algebraic closure} $\hat F$ of the free group $F$ of rank $m$.  Concordance over $\hat F$ is essentially the usual notion of link concordance. Since our understanding of $\hat F$ is very limited, it is difficult to obtain usable invariants this way (although one can map $\hat F$ to $\Z^m$ or the general lower central series quotients $F/F_q$ of $F$ (see \cite{Fr}, \cite{CK}, \cite{L}, \cite{Sm}) and use representations of these groups). For boundary links, or, more precisely $F$-links, the appropriate group is $G=F$. Since $\RR_k (F)$ is just the Cartesian product of $m$ copies of the unitary group, this gives a  workable invariant. On the other hand  concordance over $F$ is essentially what is usually called  {\em boundary concordance}, which is ostensibly stronger than ordinary concordance, It remains a major problem to decide whether these relations are really different. $\hat F$-concordance for higher dimensional links has been computed by LeDimet \cite{Le} in terms of $\GG$ groups of $\Z\hat F$ and homotopy groups of the {\em topological localization} of the wedge of circles. However, at this point we know very little about  these groups. $F$-concordance, for higher-dimensional links, has been classified by Cappell-Shaneson \cite{CS} in terms of $\GG$-groups of $\Z F$, and by Ko \cite{K} and Mio \cite{M} in terms of Seifert matrices. Recently D. Sheiham \cite{S} has reinterpreted the Ko-Mio classification and, in particular, the torsion-free part (the $F$-concordance classes of higher dimensional $F$-links form an abelian group) is measured by signatures of representations of a certain {\em quiver}. It is an obvious question, since the $\rho$-invariant is essentially a signature invariant itself,  to ask whether the $\rho$-invariant captures all of the Sheiham signatures. The main question we will address here is whether the $\rho$-invariant might be just an ordinary concordance invariant. If the answers to both questions were YES then we would conclude that concordance and boundary concordance are the same, modulo torsion, in higher dimensions. 

In this note we will apply a result of P. Vogel \cite{V} to show that the answer to the second question is YES, i.e. concordant boundary links have equivalent $\rho$-invariants.

\section{Preliminaries}
We recall some notions, working always in the smooth category. An $m$-component $n$-link $L$ is an oriented submanifold of $S^{n+2}$ which is homeomorphic to $m$ disjoint copies of $S^n$. An $m$-component $n$-disk link $D$ is a proper oriented submanifold of $D^{n+2}$ which is homeomorphic to $m$ disjoint copies of $D^n$ and such that $D\cap S^{n+1}=\bd D$ is the standard trivial $(n-1)$-link $T^{n-1}$. Thef closure of a disk link is the link obtained by adjoining to its boundary the standard trivial disk-link.

Two $m$-component $n$-links $L_0$ and $L_1$ are {\em concordant} if there is a proper submanifold $C\sub I\times S^{n+2}$, homeomorphic to $m$ copies of $I\times S^n$, such that $C\cap (t\times S^{n+2})=L_t ,\ t=0, 1$. Two $m$ component $n$-disk links are concordant if there is a proper submanifold $C$ of $I\times D^{n+2}$ such that $C\cap t\times D^{n+2}=D_t$ and $C\cap I\times S^{n+1}=I\times T^{n-1}$. We say a link, or disk link, is {\em slice} if it is concordant to the trivial link, or disk link.

The concordance classes of $m$-component $n$-disk links form a group $C(m,n)$ under a form of stacking (see \cite{Le}) This group is abelian if $n>1$ or $m=1$, but  is non-abelian if $m>1$ (see \cite{Le} again). There is a group action on $C(m,n)$ whose orbits correspond, under the closure operation, to concordance classes of links (see \cite{Le} for $n>1$ and Habegger-Lin \cite{HL} for $n=1$). Note that a disk link is slice if and only if its closure is slice. 

An $m$-component $n$-link $L$ is a {\em boundary link} if  $L$ bounds an oriented submanifold $V$ of $S^{n+2}$ consisting of $m$ components, each bounded by one of the components of $L$. $V$ is called a {\em Seifert surface} for $L$. This is equivalent to asking that there exist an epimorphism $\phi :\pp (S^3 -L)\to F$, where $F$ is the free group with a basis denoted $t_1 ,\ldots ,t_m$, such that there exist meridian elements $\mu_i\in\pp\s3$ with $\phi (\mu_i )=t_i$ \cite{G}.
$\phi$ is called an {\em $F$-structure} for $L$ and we call the pair $(L, \phi )$  an {\em $n-F$-link}. Two boundary links $L_0$ and $L_1$ are {\em boundary concordant} if there is a concordance $C$ between $L_0$ and $L_1$ such that, for some Seifert surfaces $V_t$ for $L_t$ there is a submanifold $W$ of $I\times S^{n+2}$ consisting of $m$-components, each bounded by a component of $C\cup 0\times V_0\cup 1\times V_1$. Two $F$-links $(L_t ,\phi_t )$ are $F$-concordant if there exists a concordance $C$ and an epimorphism $\Phi :\pp\sc\to F$ which restricts to $\phi_t$ on the $t$-th level. If $(L_t ,\phi_t )$ are $F$-concordant then $L_t$ are boundary concordant. The set $CF(m,n)$ of $F$-concordance classes of $m$-component $n-F$-links form an abelian group, under a suitable notion of connected sum, for $n>1$ or $m=1$. It is easy to see that $FC(1,n)=C(1,n)$, which is the usual knot concordance group. There is a group action by a certain subgroup of the group of automorphisms of $F$ on $CF(m,n)$ whose orbits correspond to boundary concordance classes of boundary links (see \cite{K}). Note that an $F$-link is $F$-slice if and only if its underlying boundary link is boundary slice.

We now point out that the natural forgetful map which associates to a boundary link  its underlying link can be lifted to a homomorphism from $F$-concordance classes of $F$-links to concordance classes of disk links.

\begin{proposition}   There is a homomorphism $\Psi_{m,n}:CF(m,n)\to C(m,n)$, for $n>1$, so that the closure of $\Psi_{m,n}(L,\phi )$ is $L$, up to concordance.
\end{proposition}

\begin{proof} Since $CF(m,2)=0$ (see \cite{K}), we may assume $n\ge 3$.

Suppose $\L =(L, \phi )$ is an $F$-link such that $\phi$ is an isomorphism---it is shown in \cite{K} that every $F$-link is $F$-concordant to such a link. We can choose  arcs $\g_i$ in $S^3$ from a base-point $*$ to each of the components of $L$, which are disjoint from each other and $L$ except at the end points, so that the meridian elements of $\pp\s3$ whose stems are the $\g_i$ are the unique elements mapped to $t_i$ by $\phi$. We can choose a regular neighborhood $R$ of $\cup_i \g_i$ which intersects each  component of $L$ in a small arc. Then $(\overline{S^3 -R},  L\cap\overline{S^3 -R})$ forms a disk link which we want to define to be a representative of $\Psi_{m,n}(\L)$. To see that this is well-defined consider $F$-cobordant $F$-links $(L_t ,\phi_t )$, where $\phi_t$ are isomorphisms, with corresponding arcs $\g^t_i$. Then we can choose an $F$-concordance $(C,\Phi )$ so that $\Phi$ is an isomorphism (see e.g. \cite{K}). It then follows that $\g^0_i$ is homotopic to $\g^1_i$ in $\I\times S^3$,where one end of the homotopy slides along $I\times *$,the other end stays in $C$ and the interior lies in $\sc$.  Since $n\ge 3$ we can perturb these homotopies to construct $m$ disjoint imbedded rectangles in $I\times S^3$ whose  sides are $\g^t_i, I\times\*$ and arcs on $C$ and otherwise is disjoint from $C$. These can be used to imbed $I\times R$ in $I\times S^3$ so that $(\overline{I\times S^3 -I\times R} , C\cap (\overline{I\times S^3 -I\times R}))$ is a concordance between $(\overline{S^3 -R}. L_0\cap (\overline{S^3 -R}))$ and $(\overline{S^3 -R}, L_1\cap (\overline{S^3 -R}))$.

Thus $\Psi_{m,n}$ is well-defined. It is easy to see that $\Psi_{m,n}$ is a homomorphism.

\end{proof}

When $n=1, CF(m,1)$ is not a group and $C(m,1)$ is a non-abelian group (unless $m=1$) and there does not seem to be any useful way to define a map $CF(m,1)\to C(m,1)$.

\section{The $\rho$-invariant.}

In \cite{APS} Atiyah-Patodi-Singer define a real-valued invariant $\rho (M, \th )$ associated to an oriented closed connected odd-dimensional manifold $M$ and a unitary representation $\th$ of $\pp (M)$. The signature theorem says that if $M$ bounds an oriented manifold $V$ and $\th$ extends to a representation $\Th$ of $\pp (V)$ then $\rho (M, \th )=k\ss  (V)-\ss (V, \Th )$, where $k$ is the dimension of the representation, $\ss (V)$ is the signature of the usual intersection pairing of $V$ and $\ss (V, \Th )$ is the signature of the intersection pairing on the twisted homology of $V$ with coefficients defined by $\Th$.

If $\L =(L, \phi )$ is an $n-F$-link we consider the closed $(n+2)$-manifold $M_L$ obtained by doing surgery on $S^{n+2}$ along the components of $L$ with the unique normal framing if $n>1$, or with the $0$-framing if $n=1$. Then $\phi$ induces an epimorphism, which we also call $\phi :\pp (M_L )\to F$. For any unitary representation $\a$ of $F$ we obtain a representation $\a\circ\phi$ of $\pp (M_L )$. If $n=2q-1$ is odd then we can take the $\rho$-invariant $\rho (M_L , \a\circ\phi )$. This defines a function
$$\rho_{\L}:\RR_k (F)\to\R $$
where $\RR_k (F)$ is the space of $k$-dimensional unitary representations of $F$, and is obviously identified with the Cartesian product of $m$ copies of the unitary group. It is shown in \cite{L} that $\rho_L$ is piecewise continuous. In fact $\rho_L$ is continuous on the complement of $\D_{\L}$, the real subvariety of $\RR_k (F)$ consisting of all $\a$ such that $H_q (M_L ; \C^k_{\a})\not= 0$, where $\C^k_{\a}$ is the local coefficient system defined by $\a\circ\phi$. We call this the {\em discriminant} of $\rho_{\L}$. It is shown in \cite[Theorem 2.1]{L} that by introducing a simple correction term, which depends only on the signatures of the component knots of $\L$, we can obtain an integer-valued function which is, therefore, piecewise constant on $\RR_k (F)-\D_{\L}$.

\section{The discriminant set of $\rho_{\L}$}

When $m=1$ and $k=1$, which is sufficient to consider for $m=1$ since every representation is a sum of $1$-dimensional representations, we have $\RR_1 (\Z )=S^1$ and $\D_{\L}$ is just the set of roots of the Alexander polynomial which lie on $S^1$. We will give an analogous description of $\dl$ for  $F$-links. 

Let $\cf$ denote the complex group algebra of $F$. A $k$-dimensional representation $\a$ of $F$ creates a right $\cf$-module which we denote by $\cka$.

\begin{lemma}\lbl{lem.hom}
$H_q (M_L ; \C^k_{\a})=0$ if and only if 
$$\cka\otimes_{\cf} H_q (\ti M_L )=\cka\otimes_{\cf} H_{q-1}(\ti M_L )=0$$
 where $\ti M_L$ denotes the $F$-cover of $M_L$ defined by $\phi$.
\end{lemma}
\begin{proof} Since $\cf$ is a free ideal ring (see \cite{C}) the usual universal coefficient theorem applies. If ${}\bf C$ is a free chain complex of left $\cf$-modules and $B$ is a right $\cf$-module, then we have a short exact sequence:

$$0\to B\otimes_{\cf}H_q ({\bf C})\to H_q (B\otimes_{\cf}{\bf C})\to\tor_1^{\cf}(B, H_{q-1}({\bf C}))\to 0$$
\begin{claim} $\tor_1^{\cf}(\cka , H_{q-1}(\ti M_L ))\iso\cka\otimes_{\cf}H_{q-1}(\ti M_L )$. 
\end{claim}
This will  complete the proof of the Lemma.

To prove the Claim we use the fact, which will be proved in Lemma \ref{lem.seif}, that $H_{q-1}(\ti M_L )$ has a free resolution $0\to M_1\to M_2\to H_{q-1}(\ti M_L )\to 0$, where $M_1$ and $M_2$ have the same rank. Tensoring this short exact sequence with $\cka$ shows that $\tor_1^{\cf}(\cka , H_{q-1}(\ti M_L ))$ and $\cka\otimes_{\cf}H_{q-1}(\ti M_L )$ are, respectively, the kernel and cokernel of a map between finite-dimensional vector spaces, and so have the same dimension (over $\C$).
\end{proof}

Let $V$ be a Seifert surface for the $F$-link $\L$. We have Seifert pairings $\a_i :H_i (V)\otimes H_{n+1-i}(V)\to\Z$ for $1\le i\le n$ and, if we choose bases for the torsion-free quotients of $H_i (V)$, we have Seifert matrices $A_i$ which represent the Seifert pairings $\a_i$ (see \cite{L3} or \cite{K} for the definitions).  We will show that the discriminant locus of $\L$  is defined by a system of equations derived from $A_q$ and $A_{q-1}$. 

Using the components of $V$ we see that each $A_i$ has a square block decomposition into $m^2$ blocks and, by Poincare duality, each of these blocks is square (see \cite{K} for details). Let $T_i$ be a square matrix with the same block configuration as $A_i$ where all the non-diagonal blocks are zero while the $j$-th diagonal block is $t_j I$ ($I$ is the identity matrix of the appropriate size). The following lemma is well-known for {\em simple} links---see e.g. \cite{G}. For any matrix $A$ let $A^T$ denote the transpose matrix.

\begin{lemma}\lbl{lem.seif} For $1\le i\le n$
 \begin{enumerate}
 \item $T_i A_i -(-1)^i A_{n+1-i}^T$ is a non-degenerate presentation matrix for $H_i (\ti M_L )$.
 \item $A_i -(-1)^i A_{n+1-i}^T$ represents the intersection pairing of $H_i (V )$ with $H_{n+1-i} (V)$.
 \end{enumerate}
 \end{lemma}

\begin{proof} This is proved by a standard argument (see e.g. \cite{G}, and \cite{L3} for the knot case) using a Mayer-Vietoris argument to obtain the exact sequence:
$$\cdots\to\cf\otimes H_i (V)\xrightarrow{d_i}\cf\otimes H_i (Y)\to H_i (\ti M_L )\to\cf\otimes H_{i-1}(V)\xrightarrow{d_{i-1}}
\cdots $$
where $Y$ is the complement of $V$. Here $d_i$ is a homomorphism of free $\cf$-modules with matrix representative $T_i A_i -(-1)^i A_{n+1-i}^T$ if $1\le i\le n$. To prove (1) we need to 
 show that $d_i$ is a monomorphism for every $0\le i\le n$. For $i=0$ we see that $d_0$ is just  the standard identification of the augmentation ideal of $\cf$ as a free left $\cf$-module of rank $m$ with generators $t_i -1$. For $i>0$,  $d_i$ is represented by a matrix over $\cf$ which, under the augmentation homomorphism $\cf\to\C$, becomes a non-singular matrix by (2). To see that this implies that $d_i$ is a monomorphism, we need to show that any square matrix $B$ over $\cf$ whose augmentation is non-singular over $\C$ is non-degenerate over $\cf$.  WLOG we can assume that the augmentation of $B$ is the identity matrix $I$. Thus $B=I-C$, where $C$ is a matrix with entries in the augmentation ideal $\I\sub\cf$. This implies that $B$ is invertible over $\cf/\I^n$ for any $n\ge 0$, $B\i =\sum_{i=0}^{n-1}C^i$. Thus if $Bv=0$, for any column vector $v$, then the entries of $v$ must lie in $\I^n$, for every $n$. But $\cap_n \I^n =0$.

The second statement of the lemma follows by a standard argument (see \cite{L3}).
\end{proof}

It is clear that if $P$ is a presentation matrix for a left $\cf$-module $M$ then we can write down a presentation matrix $\bar P$ for the vector space $\cka\otimes_{\cf}M$ in the following way. Let $\a (t_r )=\left( z_{ij}^r\right)$. Now replace every entry of the matrix $P$ by a $k\times k$-matrix by substituting the matrix $\left( z_{ij}^r\right)$ for every occurrence of $t_r$ and any scalar term $c$ by the appropriate multiple $cI$ of the $k\times k$ identity matrix. If $P$ is a square matrix then $\det \bar P=0$ if and only if $\cka\otimes_{\cf}M\not= 0$. 

Applying this to the matrices $T_i A_i -(-1)^i A_{n+1-i}^T$ for $i=q-1$ and $i=q$ results in two equations in the variables $\{ z_{ij}^r\}$ which define a complex subvariety of $\C^{mk^2}$ whose intersection with $\RR_k (F)$ defines $\dl$. 

It is not hard to show that the equations we obtain will be the same for any Seifert matrices for $\L$.

A final comment: it should be possible to derive the equations defining $\dl$ directly from the Farber trace invariants (\cite{F}, see also \cite{GL}) of $H_{q-1} (\ti M_L )$ and $H_q (\ti M_L )$.

\section{Concordance invariance of $\rho_{\L}$}

The following result is essentially proved in \cite{L}. See also \cite{Fr} for a Seifert matrix proof.
\begin{proposition}\lbl{prop.Fconc}
If $\L_0$ and $\L_1$ are $F$-concordant $F$-links then $\rho_{\L_0}=\rho_{\L_1}$ on the complement of $\D_{\L_0}\cup \D_{\L_1}$.
\end{proposition}
\begin{proof} It is shown in \cite{L} that $\rho_{\L_0}=\rho_{\L_1}$ on the complement of some special subvariety $\SS$. Suppose that $\a\in\SS$ and $\a\notin \D_{\L_0}\cup \D_{\L_1}$. Since $\SS$ is a proper subvariety of $\RR_k (F)$ there are points of $\RR_k (F)$ arbitrarily close to $\a$ not in $\SS$. So $\rho_{\L_0}=\rho_{\L_1}$ at these points and, since $\rho_{\L_0}$ and $\rho_{\L_1}$ are continuous at $\a$, it follows that $\rho_{\L_0}(\a )=\rho_{\L_1}(\a )$.
\end{proof}

The main result of this note is the following.
\begin{theorem}\label{th.conc}
\begin{enumerate}
\item Suppose that $n>1$ and $\L_0$ and $\L_1$ are $F$-links such that $\Psi_{m,n}(\L_0 )=\Psi_{m,n}(\L_1 )$. Then $\rho_{\L_0}=\rho_{\L_1}$ on the complement of $\D_{\L_0}\cup \D_{\L_1}$.
\item If $n\ge 1$ and $L$ is a boundary link which is slice, then $\rho_{\L}=0$ on the complement of $\D_{\L}$ for any $F$-structure on $L$.
\end{enumerate}
\end{theorem}
\begin{proof}
\noi (1) Let us assume for the moment that $\phi_0$ and $\phi_1$ are isomorphisms. Then the concordance between the associated disk links can be used to construct a concordance $C$ between $L_0$ and $L_1$ so that the following diagram commutes, where $i_0 , i_1$ are induced by inclusions.

\medskip
\centerline{
\divide\dgARROWLENGTH by2
$\begin{diagram}
\node[2]{\pp(S^{n+2}-L_0 )}\arrow{se,t}{i_0}\\
\node{F}\arrow{ne,t}{\phi\i_0}\arrow{se,t}{\phi\i_1}\node[2]{\pp (I\times S^{n+2}-C)}\\
\node[2]{\pp(S^{n+2}-L_1 )}\arrow{ne,b}{i_1}
\end{diagram}$}

We can now say, in the terminology of \cite{L}, that $(L_0 ,i_0 )$ and $(L_1 ,i_1 )$ are concordant $\pi$-links, where $\pi =\pp (I\times S^{n+2}-C)$. It then follows from \cite{L} that the associated $\rho$-invariants $\rho'_{\L_t}:\RR_k (\pi )\to\R,\ (t=0, 1)$ are equal outside of some special subvariety $\SS$ of $\RR_k (\pi )$. Clearly $\rho_{\L_t}=\rho'_{\L_t}\circ i^*$, where $i^* :\RR_k (\pi )\to\RR_k (F)$ is induced by $i=i_0\circ\phi\i_0 =i_1\circ\phi\i_1$. We now consider the {\em algebraic closures} $\hat F$ and $\hat\pi$ introduced in \cite{L2}.  Since $H_2 (\pi )=0$, and so $i$ is $2$-connected, it follows from \cite{L2} that $i$ induces an isomorphism $\hat F\to\hat\pi$. We now invoke a result of Vogel \cite{V} that the map $\RR_k (\hat F)\to\RR_k (F)$ induced by the natural inclusion $F\to\hat F$ is ONTO. It then follows that the map $i^* :\RR_k (\pi )\to\RR_k (F)$ induced by $i$ is onto.In fact Vogel proves the stronger fact that, if $\RR_k^o (\pi )$ is the irreducible component containing the trivial representation, then $i^* :\RR_k^o (\pi )\to\RR_k (F)$ is onto and is locally homeomorphic at the trivial representation. Thus we can conclude that $\dim\RR_k^o (\pi )=\dim\RR_k (F)$. Since $\SS'=\SS\cap\RR_k^o (\pi )$ is a proper subvariety of $\RR_k^o (\pi )$ then $\dim\SS' <\dim\RR_k (F)$ and so $i^* (\SS' )$ is a sparse subset of $\RR_k (F)$.

Now suppose $\a\in\RR_k (F)-\D_{\L_0}\cup \D_{\L_1}$. Choose $\bar\a\in\RR_k^o (\pi )$ such that $i^* (\bar\a )=\a$. If $\bar\a\notin\SS'$ then $\rho_{\L_0}(\a )=\rho'_{\L_0}(\bar\a )=\rho'_{\L_1}(\bar\a )=\rho_{\L_1}(\a )$, so let's suppose $\bar\a\in\SS'$. Since $i^* (\SS')$ is sparse in $\RR_k (F)$ there exist $\a_i\to\a$ such that $\a_i\notin i^* (\SS' )$. So if we choose $\bar\a_i\in\RR_k^o (\pi )$ over $\a_i$ then we have $\rho_{\L_0}(\a_i )=\rho'_{\L_0}(\bar\a_i )=\rho'_{\L_1}(\bar\a_i )=\rho_{\L_1}(\a_i )$. Since $\rho_{\L_0}$ and $\rho_{\L_1}$ are continuous at $\a$ we conclude that $\rho_{\L_0}(\a )=\rho_{\L_1}(\a )$.

This proves (1) with the assumption that $\phi_t$ are isomorphisms. For the general case we first recall that any $F$-link of dimension $n>2$ is $F$-concordant to a {\em simple} $F$-link (see e.g. \cite{K}), in particular  $\phi$ is an isomorphism for a simple link of dimension $>2$.
Now (1) will follow from the special case proved above, together with Proposition \ref{prop.Fconc} and the following transitivity result.
\begin{lemma}
Suppose $\L_t$ are $F$-links, $t=1,2,3$ and that $\rho_{\L_t}=\rho_{\L_{t+1}}$ outside of $\D_{\L_t}\cup \D_{\L_{t+1}}$ for $t=1,2$. Then $\rho_{\L_1}=\rho_{\L_3}$ outside of $\D_{\L_1}\cup \D_{\L_3}$.
\end{lemma}
\begin{proof} We only need to show that $\rho_{\L_1}(\a )=\rho_{\L_3}(\a )$ for $\a\in \D_{\L_2}-D_{\L_1}\cup \D_{\L_3}$. But since $\D_{\L_2}$ is a proper subvariety of $\RR_k (F)$ there are points of $\RR_k (F)$ arbitrarily close to $\a$ which are in none of the $\D_{\L_t}$ and so have the same value for every $\rho_{\L_t}$. Since $\a\notin \D_{\L_1}\cup \D_{\L_3}$ both $\rho_{\L_1}$ and $\rho_{\L_3}$ are continuous at $\a$ and so take the same value.
\end{proof}
This completes the proof of (1).

\medskip
We now prove (2). For $n>1$ this will follow from (1), so we assume $n=1$. Let $\DD\sub D^4$ be the union of slice disks for $L$. Let $\phi$ be any $F$-structure  on $L$ and set $\L =(L, \phi )$. Choose a splitting map $\mu :F\to\pp (S^3 -L)$ so that $\phi\circ\mu=$ identity, for some $F$-structure $\phi$ on $L$. Let $i:\pp (S^3 -L)\to\pp (D^4 -\DD )=\pi$ be induced by inclusion.

Since $H_2 (\pi )=H_2 (F)=0$, the maps $i, \mu ,\phi$ all induce isomorphisms of the algebraic closures of these groups. Now consider the map $\bar\phi\circ\bar i\i :\pi\to\hat F$. Recall from \cite{L2} that $\hat F$ is a direct limit of finitely presented groups $\{ G_q\}$ in such a way that the natural map $F\to\hat F$ is the direct limit of maps $\tau_q :F\to G_q$ where $\t_q$ induces an isomorphism on $H_1$ and $H_2$. It follows that we can lift $\bar\phi\circ\bar i\i$ to a map $\t :\pi\to G_q$, for some $q$ and we have the following commutative diagram:

\medskip
\centerline{
\divide\dgARROWLENGTH by2
$\begin{diagram}
\node[2]{F}\arrow{se,b}{\t_q}\arrow{ese}\\
\node{\pp (S^3 -L)}\arrow{ne,t}{\phi}\arrow{se,t}{i}\node[2]{G_q}\arrow{e}\node{\hat F}\\
\node[2]{\pi}\arrow{ne,t}{\t}\arrow{ene,b}{\bar\phi\circ\bar i\i}
\end{diagram}$}

\medskip
It now follows from \cite{L} that $\L' =(L,\t\circ i)$ is a null-concordant $G_q$-link and so the associated $\rho$-invariant $\rho_{\L'} :\RR_k (G_q )\to\R$ will vanish outside of some special subvariety $\SS$. Note that $\rho_{\L'} =\rho_{\L}\circ\t_q^*$, where $\t_q^* :\RR_k (G_q )\to\RR_k (F)$ is induced by $\t_q$. As noted above, Vogel \cite{V} shows that $\t_q^* (\RR_k^o (G_q ))=\RR_k (F)$ and, if $\SS' =\SS\cap\RR_k^o (G_q )$, then $\dim\SS' <\dim\RR_k (F)$. Thus $\t_q^* (\SS')$ is sparse in $\RR_k (F)$. 

Now to prove $\rho_{\L} (\a )=0$ for any $\a\in\RR_k (F)-\D_{\L}$ we use the same argument as in the proof of (1). Lift $\a$ to some $\bar\a\in\RR_k^o (G_q )$. If $\bar\a\notin\SS'$, the $\rho_{\L}(\a)=\rho_{\L'}(\bar\a )=0$. If $\bar\a\in\SS'$, then, since $\t_q^* (\SS' )$ is sparse in $\RR_k (F)$ we can choose $\a_i\notin\t_q^* (\SS' )$ so that $\a_i\to\a$. Lifting $\a_i$ to $\bar\a_i\in\RR_K^o (G_q )$ we have $\rho_{\L}(\a_i )=\rho_{\L'}(\bar\a_i )=0$. Since $\rho_{\L}$ is continuous at $\a$ we conclude that $\rho_{\L}(\a )=0$.

\end{proof}

\section{Problems}

We conclude with some questions.
\begin{itemize}
\item Does $\rho_{\L}$ (for all $k$) determine the $F$-concordance class of $\L$, up to torsion in the group of $F$-concordance classes, for higher dimensional $F$-links?
\item Give a direct description of the discriminant set $\D_{\L}$ using the Farber trace invariants of the link modules of $\L$.
\item Are the nilpotent representations of $F$, i.e. those that factor through $F/F_q$ for some $q$, or even those that factor through $p$-groups for some prime $p$, dense in $\RR_k (F)$? Note that this is true for $m=1$ or $k=1$. If this is true then Theorem \ref{th.conc} would follow immediately.
\end{itemize}

\ifx\undefined\bysame
	\newcommand{\bysame}{\leavevmode\hbox to3em{\hrulefill}\,}
\fi

\end{document}

%% file: abbrev.tex
\usepackage{amssymb,amsmath,amsthm,amscd,graphicx,pb-diagram,lamsarrow,pb-lams,psfrag
}
\def\draft{n}
\def\lbl#1{\label{#1}\printname{#1}}
\theoremstyle{plain}

\newtheorem{theorem}{Theorem}
\newtheorem{proposition}{Proposition}[section]
\newtheorem{lemma}[proposition]{Lemma}

\newtheorem{claim}[proposition]{Claim}

\theoremstyle{definition}

\newcommand{\noi}{\noindent}

\newcommand{\tor}{\operatorname{torsion}}

\newcommand{\tr}{\operatorname{tr}}

\def\printname#1{
	\if\draft y
		\smash{\makebox[0pt]{\hspace{-0.5in}
			\raisebox{8pt}{\tt\tiny #1}}}
	\fi
}

\def\SS{\Sigma}
\def\beq{\begin{eqnarray*}}
\def\enq{\end{eqnarray*}}
\def\ber{\left[\begin{array}{r}}
\def\berr{\left[\begin{array}{rr}}
\def\berrr{\left[\begin{array}{rrr}}
\def\berrrr{\left[\begin{array}{rrrr}}
\def\berrrrr{\left[\begin{array}{rrrrr}}
\def\berrrrrr{\left[\begin{array}{rrrrrr}}
\def\en{\end{array}\right]}
\def\a{\alpha}

\def\d{\delta}
\def\th{\theta}

\def\e{\epsilon}

\def\t{\tau}
\newcommand{\s}{\sigma}
\def\g{\gamma}
\def\l{\lambda}
\def\GG{\Gamma}
\def\DD{\Delta}

\def\Z{\mathbb Z}

\def\D{\mathsf D}

\def\C{\mathcal C}

\def\L{{\mathsf L}}
\def\R{\mathbb R}

\def\I{\mathcal I}

\def\bd{\partial}

\def\i{^{-1}}

\def\iso{\cong}

\def\ss{\smallsmile}
\def\sub{\subseteq}

\def\ti{\tilde}